\documentclass[12pt]{amsart}
\usepackage{amsbsy,amssymb,amscd,amsfonts,latexsym,amstext,delarray,
amsmath,graphicx} 
\input xypic
\usepackage[cmtip,all]{xy}

\usepackage{color}

\newtheorem{thm}{Theorem}[section]
\newtheorem{prop}[thm]{Proposition}
\newtheorem{cor}[thm]{Corollary}
\newtheorem{lem}[thm]{Lemma}

\newtheorem{defn}[thm]{Definition}

\newtheorem{rem}[thm]{Remark}

\numberwithin{equation}{section}

\def\bP{{\mathbb P}}

\def\N{{\mathbb N}}

\def\R{{\mathbb R}}

\def\cC{{\mathcal C}}
\def\cD{{\mathcal D}}
\def\cE{{\mathcal E}}
\def\cF{{\mathcal F}}

\def\cI{{\mathcal I}}
\def\cJ{{\mathcal J}}

\def\cL{{\mathcal L}}
\def\cM{{\mathcal M}}

\def\cP{{\mathcal P}}

\def\cV{{\mathcal V}}
\def\cW{{\mathcal W}}

\def\Hom{{\rm Hom}}

\newcommand\undermat[2]{%
  \makebox[0pt][l]{$\smash{\underbrace{\phantom{%
    \begin{matrix}#2\end{matrix}}}_{\text{$#1$}}}$}#2}

\title{Pareto Optimization in Categories}
\author{Matilde Marcolli}
\date{April 2022}

\address{California Institute of Technology, Pasadena \\ USA}
\email{matilde@caltech.edu}

\begin{document}
\maketitle

\begin{abstract}
We propose a model of Pareto optimization (multi-objective programming)
in the context of a categorical theory of resources. We describe how to adapt 
multi-objective swarm intelligence algorithms to this categorical formulation. 
\end{abstract}

\tableofcontents

\section{Introduction}

Pareto optimization (or multi-objective programming) 
refers to a class of problems where several simultaneous objective functions (objective
valuations), usually valued in cones inside real Euclidean spaces, need to be optimized
simultaneously. Since they are subject to constraints, the optimization cannot
be achieved simply by individually maximizing each function.  A Pareto optimal 
solution (in general non-unique) is a solution where none of the objective 
functions can be improved without worsening some of the others. More
precisely, a possible solution $S_1$ is said to Pareto dominate another solution $S_2$ 
if all the objective valuations $f_i$ satisfy $f_i(S_1)\geq f_i(S_2)$  and
for at least one of them the inequality is strict. Pareto optimal solutions
are those that are not Pareto dominated by any other. The
{\em Pareto frontier} is the set of all Pareto optimal solutions. 

\smallskip

While Pareto optimization is a very useful approach in describing optimization
problems that cannot be reduced to a single scalar function, the use of 
functions in Euclidean spaces is still often approached through a process
of  aggregative ``scalarization" that considers weighted combinations of
the different objective functions to reproduce a scalar optimization problem
for a single real valued function.

\smallskip

Our goal here is to develop a setting for Pareto optimization that 
is entirely independent of real valued functions and is formulated in
terms of the ``mathematical theory of resources" (in the sense of 
\cite{CoFrSp16} and \cite{Fr17}) in a categorical framework. In the setting we develop here,
objective functions are replaced by objective {\em functors}, and
the Pareto frontier and Pareto optimization are entirely describable
in categorical terms. The categorical setting should not be
surprising, as it is easy to see that the typical universal properties 
in category theory can be expressed in the form of optimization
problems, so that category theory is indeed a natural setting
for an astract formulation of optimization problems. 

\smallskip

We will formulate here our Pareto optimization setting in
terms of assignments of resources to a (finite) set. Thus,
the solutions of our optimization problem will be
{\em summing functors} from a category of subsets
of a finite set (which we generally think of as subsystems 
of a given system) to a symmetric monoidal category of resources. 
One can refine this setting by considering, as in
\cite{ManMar}, assignments of resources to a 
network (directed graph) using an appropriate notion
of ``network summing functors". This choice of the source
category is not necessary and can be replaced by other
categories. We choose this setting because the present
paper is part of a larger ongoing study of dynamical
assignments of resources to networks in a categorical
framework, \cite {ManMar}, \cite{Mar21}.

\smallskip

The objective functors, in turn, will be functors from
the category of summing functors to other categories
containing target goal objects. These serve the purpose
of measuring whether a given assignment of resources
to the system suffice to achieve the desired goals. 

\smallskip

We introduce the Pareto frontier in terms of an
optimization on ``convertibility of resources" in a
categorical sense, and we present a formulation of
the Pareto frontier as a category of essential preimages
of colimits of a class of diagrams.

\smallskip

In the usual setting of Pareto optimization, multi-objective programming
can be formulated in terms of a {\em swarm intelligence} algorithm. The
goal of the swarm particles in this approach is to
find solutions as close as possible to the Pareto
frontier and as diverse as possible, mapping out different regions of the Pareto frontier.
This is achieved by considering a virtual swarm of ``particles" that moves according to 
some dynamical rules across the landscape of possibilities (the configuration space).
The structure of the swarm intelligence algorithm can be summarized as follows:
\begin{enumerate}
\item the swarm is initialized by a random distribution of positions and momenta drawn with
uniform measure over configuration space;
\item each individual particle in the swarm can memorize its best solution 
up to the present time;
\item each particle in the swarm tends to search near its
best position obtained so far;
\item each particle can see the positions of the
other particles of the swarm at the same time and evaluate the best position
achieved by the swarm at that moment;
\item each individual particle tends to
move towards the best position achieved within the swarm at that time.
\end{enumerate}
This idea is formalized by update rules (in discrete time) for positions and velocities of the swarm particles,
where the velocities are updated by a rule of the form
$$ V_i(n+1) = \lambda_3\, V_i(n) + \lambda_1 \, G_1(X_i(n)-X_{i,best}(n)) + \lambda_2 G_2(X_i(n)-X_{best}(n)) $$ 
where the $G_i$ are Gaussians, the $\lambda_i$ are tunable parameters, the $X_{i,best}(n)$ is the best position of the $i$-th particle in
its previous history, that is, among the set of positions $\{ X_i(0),\ldots, X_i(n) \}$, 
while $X_{best}(n)$ is the best position among all the $N$ swarm particles $\{ X_1(n),\ldots, X_N(n) \}$  at the given time $n$.
The positions are then simply updated by the rule 
$$ X_i(n+1)= X_i(n) + V_i(n+1)\, . $$
Under good conditions (see the discussion in \cite{Mirja}, \cite{Nedja}),  
for large swarm size $N$ and sufficiently many iterations $n$, the  
positions of the swarm draw out the Pareto frontier. 

\smallskip

We show that a direct probabilistic analog of this swarm intelligence algorithm can be developed in the categorical framework. 
A single particle model identifies the Pareto frontier, but the resulting probability distribution is very spread out so it does not
provide an efficient algorithm. 

\section{Pareto frontier in categories}

\subsection{Categories of resources}

As in \cite{CoFrSp16} and \cite{Fr17}, {\em resources} are modelled by a 
symmetric monoidal category $(\cC,\circ,\otimes,{\mathbb I})$ (which we will 
also be writing ``additively" as $(\cC,\oplus, 0)$). Objects $A\in {\rm Obj}(\cC)$ 
represent resources,  the monoidal operation $A\otimes B$ represents the 
combination of resources, with the unit object ${\mathbb I}$ of the
monoidal structure representing the empty resource. 
Morphisms $f: A\to B$ in ${\rm Mor}_\cC(A,B)$ describe possible processes of
conversion of resource $A$ into resource $B$. Thus, the
convertibility of resources is expressed by the condition ${\rm Mor}_\cC(A,B)\neq \emptyset$. 

One associates to a category $(\cC,\circ,\otimes,{\mathbb I})$ of resources 
a preordered abelian semigroup
$(R,+,\succeq, 0)$ on the set $R$ of isomorphism classes of objects in ${\rm Obj}(\cC)$ with $A+B$ 
the class of $A\otimes B$ with unit $0$ given by the class of the unit object ${\mathbb I}$ 
and with $A \succeq B$ iff ${\rm Mor}_\cC(A,B)\neq \emptyset$. 
Measuring semigroups are abelian semigroups with partial ordering and with a 
semigroup homomorphism $M: (R,+)\to (S,*)$ with $M(A)\geq M(B)$ in $S$ when $A\succeq B$ in $R$.
It is shown in \cite{Fr17} that they satisfy $\rho_{A\to B} \cdot M(B) \leq M(A)$ with respect to the
maximal conversion rate 
$$  \rho_{A\to B}:= \sup \{ \frac{m}{n} \,|\, n\cdot A \succeq m \cdot B, \, \, m,n\in \N \} \, . $$

\subsection{Summing functors}

A summing functor is a {\em consistent assignment} of resources of type $\cC$
to {\em all subsystems} of a given system so that a combination of independent 
subsystems corresponds to combined resources. The notion of summing functors
was first introduced by Segal in \cite{Segal}, in the homotopy theory setting
of Gamma-spaces, in the case where $\cC$ is a category with sum and zero-object,
and was extended by Thomason in \cite{Tho82}, \cite{Tho95} to the more general case
where $\cC$ is a symmetric monoidal category. We formulate it here for finite
sets rather than for finite pointed sets as in the original setting.

\smallskip

Let $(\cC,\oplus, 0)$) be a symmetric monoidal category, written in additive notation. 
Let $S$ be a finite set and let $\cP(S)$ denote the category with objects the subsets $A\subseteq S$
and morphisms the inclusions $j: A\subseteq A'$. A functor $\Phi_S : \cP(S)\to \cC$ is a {\em summing functor} if 
$$ \Phi_S (A\cup A')=\Phi_S(A)\oplus \Phi_S (A') \ \ \  \text{ when } \ \ \  A\cap A'=\emptyset $$
and $\Phi_S (\emptyset)$ is the monoidal unit $0$ of $\cC$.

\smallskip

Let $\Sigma_\cC(S)$ be the category of summing functors $\Phi_S : \cP(S)\to \cC$, with morphisms
given by the {\em invertible} natural transformations. This category describes all the possible
assignments of resources of type $\cC$ to the subsystems of $S$, with all the possible equivalences
between such assignments. 

\smallskip

A summing functor $\Phi_S: P(S) \to \cC$ completely determined by values $\Phi_S(x):=\Phi_S(\{ x \})$ 
for $x\in S$, and the category $\Sigma_\cC(S)$ of summing functors is equivalent to the category
$\hat\cC^n$, where $n=\# S$ and where $\hat\cC$ is the category with same objects as $\cC$ 
and the invertible morphisms of $\cC$. 

\smallskip

When the finite set $S$ is replaced by a finite directed graph $G$, various notions of
``network summing functors" can be considered that generalize the setting above, We
refer the reader to \cite{ManMar} for a more detailed discussion. For the purposes of this
paper we just discuss the case of categories of summing functors $\Sigma_\cC(S)$ as
above. The generalization to networks is straightforward. 

\subsection{Objective valuation functors}

Let $S$ be a finite set as above, with $\Sigma_\cC(S)$ the category of summing functors 
for resources of type $\cC$. A {\em valuation system} $(F,X)=(F_\alpha, X_\alpha)_{\alpha\in \cI}$
consists of a finite family  $\{ \cV_\alpha \}_{\alpha \in \cI}$ of categories that 
describe possible objectives for optimization, with functors $F_\alpha: \Sigma_\cC(S) \to \cV_\alpha$
(valuations) and objects $X_\alpha \in {\rm Obj}(\cV_\alpha)$ (goals).
Valuation functors may factor through the target category of resources $\cC$, but we do not assume that
this is necessarily the case.
Valuation functors $F_\alpha: \Sigma_\cC(S) \to \cV_\alpha$ are in general {\em not} fully faithful. 

\smallskip

A summing functor $\Phi\in \Sigma_{\cC}(S)$ is {\em $F$-minorized} by another $\Psi\in \Sigma_{\cC}(S)$ if
$$ \Hom_{\cV_\alpha}(F_\alpha(\Phi),F_\alpha(\Psi))\neq \emptyset \ \ \  \forall \alpha\in \cI\, . $$
It is {\em strictly} $F$-minorized if the above holds and there exists some $\alpha\in \cI$ with 
$F_\alpha(\Phi)$ and $F_\alpha(\Psi)$ not isomorphic in $\cV_\alpha$ (hence $\Phi$ is not isomorphic
to $\Psi$). We can define $F$-majorization in a similar way.

\smallskip

The $F$-minorization condition above means that $F_\alpha(\Psi)$ is obtainable 
from $F_\alpha(\Phi)$ through an admissible ``conversion of resources" in the
category $\cV_\alpha$.

\smallskip

A summing functor $\Phi\in \Sigma_{\cC}(S)$ is {\em $(F,X)$-minorized} by another $\Psi\in \Sigma_{\cC}(S)$ if
for all $\alpha\in \cI$
$$ \begin{array}{l}\Hom_{\cV_\alpha}(F_\alpha(\Phi),F_\alpha(\Psi))\neq \emptyset \, , \\ 
 \Hom_{\cV_\alpha}(F_\alpha(\Phi),X_\alpha)\neq \emptyset  \, , \\ 
 \Hom_{\cV_\alpha}(F_\alpha(\Psi), X_\alpha)\neq \emptyset  \, , \end{array} $$
 with a strict minorization if, for some $\alpha$, $F_\alpha(\Phi)$ and $F_\alpha(\Psi)$ are
 not isomorphic. 
 This means that $F_\alpha(\Psi)$ is obtainable from $F_\alpha(\Phi)$, while both are good 
 enough to obtain the goals $X_\alpha$. 
 
 \smallskip
 
 We then define the Pareto frontier in the following way. An assignment of
 resources  $\Phi\in \Sigma_{\cC}(S)$ is on the $(F,X)$-{\em Pareto upper frontier} if
$$ \Hom_{\cV_\alpha}(F_\alpha(\Phi),X_\alpha)\neq \emptyset \ \ \  \forall \alpha\in \cI $$
but there is no $\Psi \in \Sigma_{\cC}(S)$ not isomorphic to $\Phi$ that is a strict $(F,X)$-minorization of $\Phi$.

\smallskip

The terminology ``upper frontier" is used here to indicate an optimization over valuations that lie ``above the goals". An analogous notion of lower frontier can be defined with the condition $\Hom_{\cV_\alpha}(X_\alpha,F_\alpha(\Phi))\neq \emptyset$ and $(F,X)$-majorizations.

\smallskip
\subsection{Categorical Pareto frontier}

We give here a description of the Pareto frontier as a category.

\subsubsection{Preorders and diagrams}

A preorder $\preceq$ on a set $J$ is a relation that is transitive ($x\preceq y$ and $y\preceq z$ $\Rightarrow$ $x\preceq z$) and reflexive ($x\preceq x$ for all $x$).
A preorder $(J,\preceq)$ is a directed set if $J\neq \emptyset$ and for all $x,y\in J$ there is a $z\in J$ with $x\preceq z$ and $y\preceq z$. Every finite subset $\{ x_1, \ldots, x_n \}$
of a directed set  $(J,\preceq)$ has an upper bound, that is, an element $z$ such that $x_i\preceq z$ for all $i=1,\ldots,n$. A preorder $(J,\preceq)$ is a thin category with
objects $x\in J$ and a single morphism $x\to y$ when $x\preceq y$. A directed set is a filtered thin category (all finite diagrams have a cocone). 

\smallskip

Let $\Sigma^{adm}_\cC(S)$ denote the full subcategory of $\Sigma_\cC(S)$ of $(F,X)$-admissible summing functors, with objects those
$\Phi\in {\rm Obj}(\Sigma_\cC(S))$ such that $$\Hom_{\cV_\alpha}(F_\alpha(\Phi),X_\alpha)\neq \emptyset \, , $$
for all $\alpha\in \cI$.
On ${\rm Obj}(\Sigma^{adm}_\cC(S))$ consider the preorder relations $\Psi \preceq_\alpha \Phi$
iff $$\Hom_{\cV_\alpha}(F_\alpha(\Phi),F_\alpha(\Psi))\neq \emptyset \, . $$ We write $(\cJ_\alpha, \preceq_\alpha)$
for the thin category with objects ${\rm Obj}(\Sigma^{adm}_\cC(S))$ and a single morphism $\Psi \to \Phi$
iff $\Psi \preceq_\alpha \Phi$. 

\smallskip

Let $\cD_\alpha$ denote the category with objects $Y \in {\rm Obj}(\cV_\alpha)$
that are isomorphic $Y \simeq F_\alpha(\Phi)$ in $\cV_\alpha$, for some $\Phi\in \Sigma^{adm}_\cC(S)$,
and with morphisms $\varphi=(u,v,v')$ that form a commutative diagram
$$ \xymatrix{  Y \ar[rr]^u \ar[rd]^v  & & Y' 
 \ar[ld]_{v'} \\  & X_\alpha  & } $$

\smallskip

We can then view minorizations as diagrams $D: \cJ_\alpha^{op} \to \cD_\alpha$ of the form
\begin{equation}\label{FalphaInvSys}
 \xymatrix{  F_\alpha(\Phi) \ar[rr] \ar[rd]  & & F_\alpha(\Psi)  
 \ar[ld] \\  & X_\alpha & } 
\end{equation} 
 
\smallskip
 
Thus, we are interested in considering inverse systems, namely diagrams
$D: \cJ_\alpha^{op} \to \cD_\alpha$, and their colimits in $\cD_\alpha$
$$ L_\alpha(D) := {\rm colim}_{\cJ_\alpha^{op}} D\, . $$
We consider the following class of diagrams.

\smallskip

Let ${\rm Diagr}_\alpha$ be the class of diagrams $D: \cJ_\alpha^{op} \to \cD_\alpha$ with shape
\begin{equation}\label{shape}
 \bullet \rightarrow \bullet \rightarrow \cdots \bullet \rightarrow 
\end{equation}  
either finite or infinite, with the following properties: 
\begin{enumerate}
\item if the diagram is finite of length $n$ then the object $F_\alpha(\Phi_n)$ 
in the terminal position must satisfy the condition that, for any admissible $\Psi$, 
$$ \Hom(F_\alpha(\Phi_n), F_\alpha(\Psi))\neq \emptyset \Rightarrow F_\alpha(\Phi_n)\simeq F_\alpha(\Psi)\, , $$
\item any two consecutive terms $F_\alpha(\Phi_i)\to F_\alpha(\Phi_{i+1})$ in 
the diagrams are non-isomorphic, $F_\alpha(\Phi_i)\not\simeq F_\alpha(\Phi_{i+1})$. 
\end{enumerate}

\smallskip

These are finite or infinite diagrams in $\cD_\alpha$ of the form
$$ 
 \xymatrix{  F_\alpha(\Phi_1) \ar[r] \ar[rd]  & \cdots \ar[r] \ar[d] & F_\alpha(\Phi_n) \ar[ld] & \cdots
  \\  & X_\alpha & &}
 $$
 where all the horizontal arrows are strict minorizations.
 
\smallskip 
 
 Let $\cL_\alpha =\{ L_\alpha(D) \,|\, D \in {\rm Diagr}_\alpha \}$ be the colimits of the diagrams in 
 ${\rm Diagr}_\alpha$ (when they exist in $\cD_\alpha$).  Note that, in the case of a finite diagram of
 shape $\bullet \rightarrow \cdots \rightarrow \bullet$ the colimit is isomorphic to the last term. So the
 interesting case is that of infinite diagrams.

\smallskip
\subsubsection{The Pareto frontier category}

Given a valuation functor $F_\alpha: \Sigma_\cC(S) \to \cV_\alpha$ and the class $\cL_\alpha$ of objects of $\cD_\alpha$
given by colimits of diagrams in ${\rm Diagr}_\alpha$ as above, consider the full subcategory $F_\alpha^{-1}(\cL_\alpha)$ of $\Sigma_\cC(S)$
with objects
$$ {\rm Obj}(F_\alpha^{-1}(\cL_\alpha)) =\{ \Phi\in {\rm Obj}(\Sigma_\cC^{adm}(S))\,|\, F_\alpha (\Phi) \simeq L_\alpha(D) \text{ for some } D\in {\rm Diagr}_\alpha \}\, . $$
%with objects given by those $\Phi$ that fit into a diagram in ${\rm Diagr}_\alpha$ so that the induced map $\varphi_D: F_\alpha(\Phi)\to L_\alpha(D)$ to the
%colimit is an isomorphism
%$$ {\rm Obj}(F_\alpha^{-1}(\cL_\alpha)) =\{ \Phi\in {\rm Obj}(\Sigma_\cC^{adm}(S))\,|\, 
%\varphi_D: F_\alpha (\Phi)\stackrel{\simeq}{\to} L_\alpha(D) \text{ for some } D\in {\rm Diagr}_\alpha \}\, . $$
 
\smallskip

Given the finite collection $F=\{ F_\alpha \}_{\alpha \in \cI}$ of valuation functors, and the collection
of objects $\cL=\cup_\alpha \cL_\alpha$ 
we similarly define the subcategory $\Sigma_\cC^{(F,\cL)}(S)\subset \Sigma_\cC^{adm}(S)$ 
as the full subcategory with objects given by
$$ {\rm Obj}(\Sigma_\cC^{(F,\cL)}(S)) =\cap_\alpha {\rm Obj}(F_\alpha^{-1}(\cL_\alpha)) $$

\smallskip

We can then extend our previous definition of the Pareto frontier in the following way.
We define the {\em Pareto frontier category} to be the category $\Sigma_\cC^{(F,\cL)}(S)$
obtained above. Namely,  an object $\Phi\in \Sigma_\cC^{adm}(S)$ is on the (upper) Pareto frontier
with respect to the system $(F,X)$ of valuations and goals, iff $\Phi$ is in $\Sigma_\cC^{(F,\cL)}(S)$, with $\cL_\alpha$
the colimits of diagrams in ${\rm Diagr}_\alpha$ (whenever these colimits exist in $\cD_\alpha$).

\smallskip

Note that, for all the finite diagrams, this reproduces by construction the Pareto
frontier as we described it above, since the essential preimages of the 
colimits in this case are exactly those admissible $\Phi$ that have no strict minorization, namely  
$$ \Hom(F_\alpha(\Phi), F_\alpha(\Psi))\neq \emptyset \Rightarrow F_\alpha(\Phi)\simeq F_\alpha(\Psi)\, . $$
The difference here is that we include the colimits of the infinite sequences of strict minorizations,
if these colimits exist in $\cD_\alpha$, and we describe the Pareto frontier as a category rather
than a set/class. 

\medskip
\section{Probabilistic particles}\label{probpartSec}

In this section we discuss a direct analog, in our categorical setting, of the usual swarm intelligence algorithm
for multi-objective optimization, and we show that this simple generalization does not suffice to identify the
Pareto frontier. 

\smallskip
\subsection{Probabilistic categories}

Let $\cF\cP$ be the category of finite probabilities, with objects $(X,P)$ consisting of a finite set $X$ with a probability measure $P$, 
and morphisms $S\in \Hom_{\cF\cP}((X,P),(Y,Q))$ given by stochastic $(\# Y \times \# X)$-matrices $S$, namely matrices with
$S_{yx}\geq 0$, for all $x\in X$, $y\in Y$ and $\sum_{y\in Y} S_{yx}=1$ for all $x\in X$, such that 
the probability measures are related by $Q = S\, P$.

\smallskip

As shown in \cite{Mar19}, given a category $\cC$, one can construct a probabilistic version $\cP\cC$, which can
be viewed as a wreath product $\cF\cP \wr \cC$ of the category $\cC$ with the category 
$\cF\cP$ of finite probabilities.

\smallskip

The objects of $\cP\cC$ and formal finite convex combinations
$$ \Lambda C = \sum_i \lambda_i C_i $$
with $\Lambda=(\lambda_i)$ a finite probability and $C_i \in {\rm Obj}(\cC)$.

\smallskip

Morphisms in $\Hom_{\cP\cC}(\Lambda C,\Lambda' C')$ are pairs $(S,f): \Lambda C \to \Lambda' C'$ with
$S$ a stochastic matrix with $S\Lambda =\Lambda'$, and $f=\{ f_{ab,r} \}$ 
finite collection of morphisms $f_{ab,r}: C_b \to C'_a$ with assigned probabilities 
$\mu^{ab}_r$. These probabilities satisfy $\sum_r \mu^{ab}_r =S_{ab}$.

\smallskip

Morphisms $(S,f)$ in $\cP\cC$ can be interpreted as ways of
mapping $C_b$ to $C_a'$
by randomly choosing a morphism from the set $\{ f_{ab,r} \}$,  
with probability $\mu^{ab}_r$ of choosing $f_{ab,r}$. 

\smallskip

We have the following simple characterization of isomorphic objects in probabilistic categories. 

\begin{lem}\label{isoPC}
Two objects $\Lambda X=\sum_i \lambda_{i=1}^n X_i$ and $\Lambda' X'=\sum_{j=1}^m \lambda'_j X'_j$ in a probabilistic category 
$\cP\cC$ are isomorphic if and only if $n=m$ with $X_i \simeq X'_{\sigma(i)}$ (isomorphic in $\cC$) 
for some permutation $\sigma$ and $\lambda_i=\lambda'_{\sigma(i)}$.
\end{lem}

\proof The isomorphism $\Lambda X\simeq \Lambda' X'$ means that there is an invertible morphism $(S,f): \Lambda X \to \Lambda' X'$.
In particular $S$ with $S \Lambda =\Lambda'$ must be a stochastic matrix with stochastic inverse, hence we have $n=m$ and $S$ is 
necessarily a permutation matrix. Thus, $\lambda_i=\lambda'_{\sigma(i)}$ for a permutation $\sigma$. The collection of morphisms 
$f=\{ f_{ij,r} \}$ then have probabilities $\mu_{ij,r}$ satisfying $\sum_r \mu_{ij,r}=S_{ij}$ hence they can be nonzero only for $j=\sigma(i)$,
with $\sum_r \mu_{i\sigma(i), r}=1$, and $f_{i,\sigma(i),r}: X_i \to X'_{\sigma(i)}$ an isomorphism. 
\endproof

\smallskip

For $\cC$ a small category, it is also natural to assume that an object $\sum_i \lambda_i C_i$ of $\cP\cC$ where
$C_i = C$ for all $i$ would be the same as the object $C$ with probability $\Lambda=\{ 1 \}$. However, it is better to
just require, more generally, that, whenever $C_i \simeq C$, the objects $\sum_i \lambda_i C_i$ and $C$ are isomorphic 
objects. This can be achieved by a localization of the category $\cP\cC$.

\begin{lem}\label{locPC} 
Let $\cW$ be the class of morphisms is $\cP\cC$ of the form $\varphi: \sum_i \lambda_i C_i \to C$ with 
$\varphi=(S,f)$ with an $n\times (n+k)$ stochastic matrix $S$ of the form
$$ S=\left(\begin{array}{cccc|ccc|cccc}
1& 0 & \cdots & 0 & 0 & 0 \cdots 0 & 0 & 0 & \cdots & 0 \\
0& 1 & & 0 & 0 & 0 \cdots 0 & 0 & 0 & \cdots & 0 \\
0 & 0 & \cdots & 1 & 0 & 0 \cdots 0 & 0 & 0 & \cdots & 0 \\
0 & 0 & \cdots & 0 & 1 & 1 \cdots 1 & 1 & 0 & \cdots & 0 \\
0 & 0 & \cdots & 0 & 0 & 0 \cdots 0 & 0 & 1& \cdots & 0 \\
0 & 0 & \cdots & 0 & 0 & 0 \cdots 0 & 0 & 0 & \cdots & 1 \\
\undermat{\ell}{0 & 0 & \cdots & 0} & \undermat{k}{0 & 0 \cdots 0 & 0 } & \undermat{n-\ell}{0 & \cdots & 1}  \\
\end{array}\right) $$ 

\bigskip
\noindent and where the set $f$ of morphisms consists of the identity ${\rm id}_{C_i}$ for $i=1,\ldots, \ell$
and $i=\ell+k+1,\ldots, n+k$, and of isomorphisms $f_i : C_i \to C$, for $i=\ell+1,\ell+k$, all of them occuring with probability $1$.
Then the localization $\cP\cC[\cW^{-1}]$ implements the equivalence relation described above.
\end{lem}

\smallskip

\proof The category $\cP\cC[\cW^{-1}]$ is the localization of $\cP\cC$ at $\cW$. Isomorphisms in $\cP\cC[\cW^{-1}]$ are
arbitrary compositions of the isomorphisms of Lemma~\ref{isoPC} and morphisms in $\cW$ and their formal inverses.
Thus, in this category we have isomorphisms between an object $\sum_{i=1}^m \lambda_i C_i$ where the $C_i$
for a subset $i\in I\subset \{ 1, \ldots, m \}$ of indices are all isomorphic to the same object $C$ and the object
$(\sum_{i\in I}\lambda_i) \, C + \sum_{i\in I^c} \lambda_i C_i$.
\endproof

\smallskip
\subsubsection{Probabilistic categories of functors}

Of particular interest here is the case where the category $\cC$ is a category of
functors $\cC={\rm Func}(\cD,\cD')$. In this case, when we form the probabilistic
category $\cP{\rm Func}(\cD,\cD')$, we want to interpret an object $\Lambda F=\sum_i \lambda_i F_i$ 
in $\cP{\rm Func}(\cD,\cD')$ as a functor that to an object $X$ in ${\rm Obj}(\cD)$
assigns $F_i(X) \in {\rm Obj}(\cD')$ with probability $\lambda_i$. In order to make this
heuristics precise, so that we can use it in defining swarm dynamics in categories,
we need the following simple statements.

\smallskip

\begin{lem}\label{Fprob}
Functors $F: \cC \to \cC'$ extend to functors $F: \cP \cC \to \cP \cC'$ mapping
an object $\Lambda C=\sum_i \lambda_i C_i$ of $\cP\cC$ to the object $\Lambda F(C)=\sum_i \lambda_i F(C_i)$
of $\cP\cC'$ and a morphism $(S,f)$ with $f=\{ f_{ab,r} \}$ and probabilities $\mu_{ab,r}$ to
the morphism $F(S,f)=(S,F(f))$ with $F(f)=\{ F(f_{ab,r}) \}$ with probabilities $\mu_{ab,r}$. 
\end{lem}

\smallskip

This follows directly from the definition. Moreover, there is a functor between
the probabilistic category of functors $\cP{\rm Func}(\cD,\cD')$ and the category of
functors between the probabilistic categories ${\rm Func}(\cP\cD,\cP\cD')$.

\smallskip

\begin{lem}\label{PEndPC}
There is a faithful functor $\cP{\rm Func}(\cD,\cD') \to {\rm Func}(\cP\cD,\cP\cD')$.
\end{lem}

\proof
Let $\Lambda F=\sum_i \lambda_i F_i$ be an object in the probabilistic category $\cP{\rm Func}(\cD, \cD')$. We first show that $\Lambda F$ defines a functor in ${\rm Func}(\cP\cD, \cP\cD')$.
For $\Omega X=\sum_a \omega_a X_a \in {\rm Obj}(\cP\cD)$ we take $\Lambda F( \Omega X)=\sum_{i,a} \lambda_i \omega_a F_i(X_a)=: \Lambda\Omega\, F(X)$ in
${\rm Obj}(\cP\cD')$, where $(\Lambda \Omega)_{i,a}=\lambda_i  \omega_a$. 
For $(S,f): \Omega X \to \Omega' X'$ in $\Hom_{\cP\cD}(\Omega X,\Omega' X')$, with $f=\{ f_{ab,r} \}$ with probabilities $\mu_{ab,r}$,  we set 
$\Lambda F(S,f): \Lambda\Omega\, F(X) \to \Lambda\Omega' \, F(X')$ with $\Lambda F(S,f)=(S', f')$, 
where $S'_{ij,ab}=\delta_{ij} S_{ab}$, so that $S' \, \Lambda\Omega=\Lambda\Omega'$, and $f'=\{ F_i(f_{ab,r}) \cdot \delta_{ij}  \}$ with probabilities $\mu_{ab,r}$.  
Consider then a morphism $(R,\alpha)$ in $\Hom_{\cP{\rm Func}(\cD,\cD')}(\Lambda F, \Lambda' F')$, with a stochastic matrix $R$ with $R \Lambda =\Lambda'$
and $\alpha =\{ \alpha_{ij, s} \}$ with probabilities $\nu_{ij,s}$ with $\sum_s \nu_{ij,s}=R_{ij}$, and with $\alpha_{ij,s}: F_i \to F'_j$ a collection of natural transformations between
functors in  ${\rm Func}(\cD, \cD')$. Then $(R,\alpha)$ defines a natural transformation of functors in ${\rm Func}(\cP\cD, \cP\cD')$, by taking, for each object $\Omega X\in {\rm Obj}(\cP\cD)$
assigns the morphism in $\cP\cD'$
$$ (R,\alpha)|_{\Omega X}: \Lambda F(\Omega X) \to  \Lambda' F'(\Omega X) $$
with stochastic matrix $R$ and with the collection  $\{ \alpha_{ij, s}|_{X_a} \}$ with probabilities $\nu_{ij,s}$, where $\alpha_{ij, s}|_{X_a}: F_i(X_a) \to F'_j(X_a)$
is the morphism in $\Hom_{\cD'}(F_i(X_a), F'_j(X_a))$ specified by the natural transformation $\alpha_{ij, s}$. The morphism $(R,\alpha)$ in
$\cP{\rm Func}(\cD,\cD')$ uniquely specifies this natural transformation. 
\endproof

\smallskip

We can interpret the difference between viewing an object $\Lambda F=\sum_i \lambda_i F_i$ 
with $F_i\in {\rm Obj}(\cE(\cC))$ and $\Lambda$ a probability distribution as objects of
$\cP(\cE(\cC))$ or (through the functor of Lemma~\ref{PEndPC}) as objects in $\cE(\cP(\cC))$
as, respectively, the {\em probabilistic} and {\em deterministic} interpretations of
$\Lambda F$.

\smallskip
\subsubsection{Other probabilistic conditions}

We assume in this section that the category $\cC$ or resources, where summing functors $\Phi\in \Sigma_\cC(S)$ take values, 
is a {\em small} category endowed with a {\em probability distribution} $\bP$ on the set ${\rm Obj}(\cC)$. This probability can be
seen as modeling the relative abundance or scarcity of resources. 

\smallskip

Through the identification of summing functors in $\Sigma_\cC(S)$ with objects in $\hat\cC^n$, with $n=\# S$, we then obtain an 
induced probability, which we also denote by $\bP$, on ${\rm Obj}(\Sigma_{\cC}(S))$.

\smallskip

Consider, as above, the subcategory $\Sigma^{adm}_{\cC}(S)$. The condition
$$ \bP({\rm Obj}(\Sigma^{adm}_{\cC}(S))) > 0 $$
ensures that the set $(F,X)$ of goals and valuations is not incompatible with the availability of resources of type $\cC$. 

\smallskip

Let $\cM^{adm}_{(X,F)}(\Phi)\subset {\rm Obj}(\Sigma^{adm}_{\cC}(S))$ denote the set of all strict $(F,X)$-minorizations of
$\Phi\in \Sigma^{adm}_\cC(S)$. The condition that $\Phi$ is on the Pareto frontier is then that $\cM^{adm}_{(X,F)}(\Phi)=\emptyset$.
If the measure $\bP$ has no non-empty sets of measure zero, then  $\lambda(\Phi):=\bP(\cM_{(X,F)}^{adm}(\Phi)) =0$ iff $\Phi$ is
on the Pareto frontier. 

\smallskip
\subsection{Single particle}

For the dynamics of a single particle, we initialize at time zero by drawing an object $\Phi_0$ 
from ${\rm Obj}(\Sigma^{adm}_{\cC}(S))$ uniformly at random with respect to the
probability measure $\bP$. 

\smallskip

The dynamics then proceeds by making new random steps and comparing them (``velocities" are 
here regarded as probabilistic jumps to a new position). Thus, at the first step (time $t=1$) a new draw of an
element $\Phi_1$ is performed. With probability $\lambda_0=\bP(\cM_{(X,F)}^{adm}(\Phi_0))$ this
new point improves the position with respect to $\Phi_0$, being a strict minorization of $\Phi_0$. 
If it is not (with probability $1-\lambda_0$) then one keeps the same position $\Phi_0$. This means
that the result of the first step is an object in the probabilistic category $\cP\Sigma_\cC^{adm}(S)$ of
the form
$$ (\Lambda \Phi)_1 := (1-\lambda_0) \Phi_0 + \lambda_0 \Phi_1 \, . $$
At the second step (time $t=2$), one makes another random draw $\Phi_2$. With
$\lambda_1=\bP(\cM_{(X,F)}^{adm}(\Phi_1))$, one then obtains a new object
in the probabilistic category of the form
$$ (\Lambda \Phi)_2 = (1-\lambda_0) ((1-\lambda_0) \Phi_0 + \lambda_0 \Phi_2) + \lambda_0 ((1-\lambda_1) \Phi_1 + \lambda_1 \Phi_2) $$ $$ =
(1-\lambda_0)^2 \Phi_0 + \lambda_0 (1-\lambda_1) \Phi_1 + \lambda_0 (1-(\lambda_0-\lambda_1)) \Phi_2 \, ,$$
that describes all the possible relative positions of $\Phi_0,\Phi_1,\Phi_2$ with the respective probabilities.

\smallskip

An inductive argument shows that one obtains the following behavior of this single particle
case.

\begin{prop}\label{probPhin}
After $n$ steps the outcome is an object of $\cP\Sigma_\cC^{adm}(S)$
of the form
$$ (\Lambda \Phi)_n = \sum_{k=0}^n c_n^k \Phi_k  $$ 
 with the probability $\Lambda_n=( c_n^k )_{k=0}^n$ satisfying the recursion (with $c_0^0=1$)
 \begin{equation}\label{cnkrecursion}
 \left\{ \begin{array}{ll} c_n^k = c_k^k\, (1-\lambda_k)^{n-k} & 0\leq k\leq n-1 \\[3mm] 
c_n^n = \sum_{k=0}^{n-1} \lambda_k\, (1-\lambda_k)^{n-1-k}\, c_k^k &  \end{array} \right. \, . 
\end{equation}
\end{prop} 

\proof We have $c_0^0=1$ with $c^0_1=1-\lambda_0$ and $c_1^1=\lambda_0$ as above.
At the $(n+1)$-st step we are comparing the new draw $\Phi_{n+1}$ with each of the previous $\Phi_k$:
it will be better than $\Phi_k$ (a strict minorization) with probability $\lambda_k$ and not better
with probability $1-\lambda_k$. When applied to the previous $(\Lambda \Phi)_n$ we then
obtain a new $(\Lambda \Phi)_{n+1}$ of the form
$$ c^0_n ((1-\lambda_0)\Phi_0 +\lambda_0 \Phi_{n+1}) + \cdots +
c_n^n ((1-\lambda_n) \Phi_n +\lambda_n \Phi_{n+1}) \, ,  $$
which gives $c_{n+1}^0=(1-\lambda_0)^{n+1}$, $c_{n+1}^1=\lambda_0 (1-\lambda_0)^{n}$,
$c_{n+1}^2=c_2^2 (1-\lambda_2)^{n-1}$, $\ldots$, $c_{n+1}^{n}=c_{n}^{n} (1-\lambda_n)$, and
$$ c_{n+1}^{n+1}= c_{n}^0 \lambda_0 + c_{n}^1 \lambda_1+\cdots + c_n^n \lambda_n \, . $$
The first relations give $c_{n+1}^k =c_k^k (1-\lambda_k)^{n+1-k}$, while the last one combined with this 
gives the second recursive relation of the statement.
The recursion directly implies that the normalization $\sum_k c_n^k=1$ holds.
\endproof

\smallskip

We have the following easy reformulation of the recursion \eqref{cnkrecursion}.

\begin{cor}\label{Scnk}
The recursive relation for the probabilities $c_n=(c_n^k)_{k=1}^n$ is implemented by $c_{n+1}=S_n \, c_n$,  with the $(n+1)\times n$ stochastic matrix
$$ S_n =\left( \begin{array}{ccccc}
1-\lambda_0 & 0 & 0 & \cdots & 0 \\
0 & 1-\lambda_1 & 0 & \cdots & 0 \\
0 & 0 & 1-\lambda_2 & \cdots & 0 \\
0 & 0 & 0 & \cdots & 1-\lambda_n \\
\lambda_0 & \lambda_1 & \lambda_2 & \cdots & \lambda_n 
\end{array} \right)\, . $$
\end{cor}

\smallskip

All the coefficients $c_n^k$ are polynomials in the $\lambda_i$ so they depend on $\Phi_0, \ldots, \Phi_n$.
The coefficient $c_n^k$ is the probability of having $\Phi_k$ as the ``best position" of the particle 
during the first $n$ steps. More precisely, the coefficient $c_n^k$ measures the probability that, among the
first draws $\{ \Phi_0, \ldots, \Phi_n \}$ there is a subsequence of $k$ strict minorizations ending with $\Phi_k$.

\smallskip

Simple numerical examples with different choices of a sequence $\lambda_0\geq \lambda_1\geq\cdots \geq \lambda_n \cdots$ 
show that this probability distribution $\Lambda_n=( c_n^k )_{k=0}^n$ can be 
very spread out: it peaks somewhere in $k$, but not always at the end term and can be 
very non-concentrated.

\smallskip

\begin{lem}\label{pinestimate}
For $\lambda_0\geq \lambda_1\geq \cdots \geq \lambda_n \cdots$,
the coefficients $c_n^n$ satisfy the estimate
\begin{equation}\label{estpin}
 c_k^k \, (1-\lambda_0)^{n-k} \leq c_n^n \leq c_k^k
\end{equation} 
for all $k=0,\ldots, n-1$.
\end{lem}

\proof By the recursion  \eqref{cnkrecursion} we have
$$ c_n^n =\sum_{k=0}^{n-2} c_k^k (1-\lambda_k)^{n-k}\lambda_k + c_{n-1}^{n-1} \lambda_{n-1} 
 \leq ( \sum_{k=0}^{n-2} c_k^k (1-\lambda_k)^{n-1-k}\lambda_k ) (1-\lambda_{n-1}) + c_{n-1}^{n-1} \lambda_{n-1}= c_{n-1}^{n-1} $$
while similarly we also get $c_{n-1}^{n-1} (1-\lambda_0)\leq c_n^n$. Iterating these
estimates we get \eqref{estpin}. 
\endproof

In particular we have the rough estimate $c_n^n \geq (1-\lambda_0)^n \geq (1-\lambda)^n$, where
$\lambda_0=\bP( \cM_{(X,F)}^{adm}(\Phi_0))$ and $\lambda=\bP({\rm Obj}(\Sigma^{adm}_{\cC}(S)))$. 

\smallskip

In a similar way, we can compute other relevant probabilities. For example we have the following.

\begin{lem}\label{onlyseq}
The probability that, among the first draws $\{ \Phi_0, \ldots, \Phi_n \}$, the 
subsequence $\{ \Phi_0, \Phi_{\ell_1}, \ldots, \Phi_{\ell_k} \}$ with $0\leq k\leq n$ is the maximal subsequence consisting 
of strict minorizations is given by
$$ \pi_{\ell_1,\ldots, \ell_k} :=(1-\lambda_0)^{\ell_1-1} \lambda_0 (1-\lambda_{\ell_1})^{\ell_2-\ell_1-1} \lambda_{\ell_2} \cdots  (1-\lambda_{\ell_{k-1}})^{\ell_k -\ell_{k-1}-1} \lambda_{\ell_{k-1}} 
(1-\lambda_{\ell_k})^{n-\ell_k} \, . $$
\end{lem} 

\proof This follows directly from the recursive construction of the $(\Lambda \Phi)_n$ discussed above.
\endproof

\smallskip

The expressions $\pi_{\ell_1,\ldots, \ell_k}(\lambda_0, \lambda_{\ell_1},\ldots, \lambda_{\ell_k})$ of Lemma~\ref{onlyseq} 
describe the probability that in a 
random draw of a sequence $\{ \Phi_0, \ldots, \Phi_n \}$ of objects in $\Sigma^{adm}_{\cC}(S)$ we can
form a longest chain of strict minorizations $$\{ \Phi_0, \Phi_{\ell_1}, \ldots, \Phi_{\ell_k} \}\, , $$ hence diagrams
in ${\rm Diagr}_\alpha$ of the form
$$ \xymatrix{ F_\alpha(\Phi_0) \ar[r] \ar[rrd]  & F_\alpha(\Phi_{\ell_1}) \ar[r] \ar[rd] & \cdots \ar[r] & F_\alpha(\Phi_{\ell_k}) \ar[ld] & \cdots \\
& & X_\alpha & & } $$
where the sequence $\lambda_r=\bP( \cM_{(X,F)}^{adm}(\Phi_r) )$ in this case must satisfy $\lambda_0\geq \lambda_{\ell_1} \geq \cdots \geq \lambda_{\ell_k}$, since in the case
of successive minorizations $\cM_{(X,F)}^{adm}(\Phi_{\ell_r})\subseteq \cM_{(X,F)}^{adm}(\Phi_{\ell_{r-1}})$. 
 
\smallskip
\subsubsection{Single particle and colimits}

Consider the case where a sequence $\{ \Phi_k \}_{k=0}^\infty$ of objects in $\Sigma^{adm}_{\cC}(S)$  forms an infinite system $D$ of strict minorizations in ${\rm Diagr}_\alpha$ of the form
\begin{equation}\label{infiniteD}
 F_\alpha(\Phi_0) \stackrel{\varphi_0}{\rightarrow} F_\alpha(\Phi_1) \stackrel{\varphi_1}{\rightarrow} \cdots \rightarrow F_\alpha(\Phi_n)\stackrel{\varphi_n}{\rightarrow} \cdots  
\end{equation}
that has a colimit  $L_\alpha(D)$ in $\cD_\alpha$. 

\smallskip

\begin{prop}\label{colimPcolim}
\begin{enumerate}
\item The system \eqref{infiniteD} induces a system in $\cP\cD_\alpha$ of the form
\begin{equation}\label{infinitePD}
\cdots \rightarrow \sum_{k=0}^n c_n^k \, F_\alpha(\Phi_k) \stackrel{(S_n,\underline{\varphi}_n)}{\rightarrow}  \sum_{k=0}^{n+1} c_{n+1}^k \, F_\alpha(\Phi_k) \rightarrow \cdots 
\end{equation}
with the $S_n$ are as in Corollary~\ref{Scnk} and the set $\underline{\varphi}_n$ consists of the maps ${\rm id}_{F_\alpha(\Phi_k)}$ for $k=0,\ldots, n$ with
probability $1-\lambda_k$ and $\varphi_{k,n+1}=\varphi_n \circ \cdots \circ \varphi_k$ with probability $\lambda_k$. 
\item A collection of $M$ cocones in $\cD_\alpha$ given by commutative diagrams
\begin{equation}\label{diagrYr}
 \xymatrix{ F_\alpha(\Phi_0) \ar[r]^{\varphi_0} \ar[rrd]^{f_{0,r}}  & F_\alpha(\Phi_{\ell_1}) \ar[r]^{\varphi_1} \ar[rd]^{f_{1,r}} & \cdots \ar[r] & F_\alpha(\Phi_{\ell_k}) \ar[r]^{\varphi_k} \ar[ld]^{f_{n,r}} & \cdots \\
& & Y_{\alpha,r} & & } 
\end{equation}
for $r=1,\ldots, M$, together with a sequence of $M\times n$ stochastic matrices $\tilde S_n$ with columns equal to the uniform distribution $1/M$ on the set $r\in \{ 1, \ldots, M \}$, 
induces a cocone in $\cP\cD_\alpha$ with commutative diagrams
\begin{equation}\label{diagrYM}
  \xymatrix{ \cdots \ar[r]  & \sum_{k=0}^n c_n^k \, F_\alpha(\Phi_k)  \ar[r]^{(S_n,\underline{\varphi}_n)} \ar[rd]^{(\tilde S_n, \underline{f}_n)} &  \sum_{k=0}^{n+1} c_{n+1}^k \, F_\alpha(\Phi_k) \ar[r] \ar[d]^{(\tilde S_{n+1},  \underline{f}_{n+1})} & \cdots \\
& & \tilde\Lambda Y_\alpha & & } 
\end{equation}
where $\tilde\Lambda Y_\alpha=\sum_{r=1}^M \tilde\lambda_r Y_{\alpha,r}$ with $\tilde\Lambda=(\tilde\lambda_j=1/M)$ the uniform probability
distribution $1/M$. The set of morphisms 
$\underline{f}_n$ consists of the $f_{k,r}: F_\alpha(\Phi_k)\to Y_{\alpha,r}$, for $k=0,\ldots, n$, each occurring with probability $1/M$. 
\end{enumerate}
\end{prop}

\proof (1) This follows from the fact that the maps $\varphi_{k,n+1}=\varphi_n \circ \cdots \circ \varphi_k$ satisfy the properties of a directed system
$\varphi_{n+2, m+1}\circ \varphi_{k,n+1}=\varphi_{k,m+1}$ and $\varphi_{k,k}={\rm id}_{F_\alpha(\Phi_k)}$.

(2) In order to obtain commutative diagrams the sequence of stochastic matrices $\tilde S_n$ must satisfy the recursive condition
$\tilde S_{n+1} \cdot S_n = \tilde S_n$, or equivalently $(1-\lambda_k) (\tilde S_{n+1})_{rk}  + \lambda_k (\tilde S_{n+1})_{r\, n+1}= (\tilde S_n)_{rk}$,
for $k=1,\ldots, n$. This recursive condition $\tilde S_{n+1} \cdot S_n = \tilde S_n$ is solved by the stochastic matrices $\tilde S_n$ with
columns given by the uniform distribution $1/M$ on $r\in \{ 1, \ldots, M \}$. Thus, we obtain commutative diagrams 
$$  \xymatrix{ \cdots \ar[r]  & \sum_{k=0}^n c_n^k \, F_\alpha(\Phi_k)  \ar[r]^{(S_n,\underline{\varphi}_n)} \ar[rd]^{(\tilde S_n, \underline{f}_n)} &  \sum_{k=0}^{n+1} c_{n+1}^k \, F_\alpha(\Phi_k) \ar[r] \ar[d]^{(\tilde S_{n+1},  \underline{f}_{n+1})} & \cdots \\
& & \tilde\Lambda Y_\alpha & & } $$
where $\tilde\Lambda Y_\alpha=\sum_r \tilde\lambda_r Y_{\alpha,r}$, with the probability measure $\tilde \Lambda$ determined by $\tilde S_n c_n =\tilde\Lambda$, where the
left-hand-side is independent of $n$ by the conditions $\tilde S_{n+1} \cdot S_n = \tilde S_n$  and $S_n c_n =c_{n+1}$. For $\tilde S_n$ as above this gives that $\tilde \Lambda$ is the uniform
distribution $1/M$. This diagram defines a cocone in $\cP\cD_\alpha$.  
\endproof

\smallskip

\begin{cor}\label{colimPD}
Consider a diagram \eqref{infiniteD} with colimit  $L_\alpha(D)$ in $\cD_\alpha$, and the class of cocone diagrams in $\cP\cD_\alpha$ of the form \eqref{diagrYM}
obtained as in Proposition~\ref{colimPcolim}. This class of cocones, seen in the localization $\cP\cD_\alpha[\cW^{-1}]$, has colimit isomorphic to $L_\alpha(D)$.
\end{cor}

\proof For each diagram we have canonical maps from the system to the colimit and from the colimit to the tip of the cocone
$$ \xymatrix{ F_\alpha(\Phi_0) \ar[r]^{\varphi_0} \ar[rrd] \ar[rrdd] & F_\alpha(\Phi_{\ell_1}) \ar[r]^{\varphi_1} \ar[rd] \ar[rdd] & \cdots \ar[r] & F_\alpha(\Phi_{\ell_k}) \ar[r]^{\varphi_k} \ar[ld] \ar[ldd] & \cdots \\
& & L_\alpha(D) \ar[d] & & \\
& & Y_{\alpha, r} & &
} $$
As in Proposition~\ref{colimPcolim} this gives induced diagrams in $\cP\cD_\alpha$ with $\tilde\Lambda$ the uniform probability distribution,
$$ \xymatrix{ \cdots \ar[r]  & \sum_{k=0}^n c_n^k \, F_\alpha(\Phi_k)  \ar[rr]^{(S_n,\underline{\varphi}_n)} \ar[rd] \ar[rdd] & &  \sum_{k=0}^{n+1} c_{n+1}^k \, F_\alpha(\Phi_k) \ar[r] \ar[ld] \ar[ldd]
& \cdots \\
& & \tilde\Lambda L_\alpha(D) \ar[d] & & & \\
& & \tilde\Lambda Y_\alpha & & &} $$
hence we can identify $\tilde\Lambda L_\alpha(D)$ with the colimit of this class of cocones.
By viewing these induced diagram in the localization $\cP\cD_\alpha [\cW^{-1}]$ we obtain an isomorphism $\tilde\Lambda L_\alpha(D) \simeq L_\alpha(D)$. 
\endproof

\smallskip

In particular, consider the case of finite diagrams in $D$ of strict minorizations in ${\rm Diagr}_\alpha$, with $\Phi_n$ on the Pareto frontier,
\begin{equation}\label{finDn}
 \xymatrix{ F_\alpha(\Phi_0) \ar[r] \ar[rrd]  & F_\alpha(\Phi_1) \ar[r] \ar[rd] & \cdots \ar[r] & F_\alpha(\Phi_n) \ar[ld]  \\
& & X_\alpha & & } 
\end{equation}
where we have $\lambda_0\geq \lambda_1\geq \cdots \geq \lambda_{n-1}$ and $\lambda_n=0$, since $\Phi_n$ is
on the Pareto frontier, hence $\cM_{(X,F)}^{adm}(\Phi_n)=\emptyset$. We obtain in this case, by the same argument above,
that the colimit of the induced finite diagrams of the $\sum_k c_m^k F_\alpha(\Phi_k)$ with $0\leq m \leq n$, is
given by the object $\tilde\Lambda F_\alpha(\Phi_n)$ with $\tilde\Lambda$ the uniform distribution. This is isomorphic to $F_\alpha(\Phi_n)$
in the localization $\cP\cD_\alpha[\cW^{-1}]$.

\smallskip

This shows that our probabilistic single particle model still computes the same colimit over a sequence of minorizations in ${\rm Diagr}_\alpha$, so it does identify
the Pareto frontier, both in the case of finite and of infinite chains of minorizations. 

\smallskip

However, because of the fact that the probability distribution $c_n=\{ c_n^k \}_{k=1}^n$ tends to be very spread out, the single particle model does not provide an
efficient computational method. Moreover, we have not established yet a method of construction of successive approximations to the Pareto frontier using
sequences of random draws $\{ \Phi_n \}_{n\geq 0}$ and their associated probabilistic objects $\sum_k c_n^k \Phi_k$. We will address these in the next section.

\medskip
\section{Swarms in categories}\label{swarmSec}

In order to address approximation, we need to assume additional structure on the target categories $\cV_\alpha$ of the
valuation functors. In particular, we introduce a {\em scale} parameter in the $\cV_\alpha$ so that we can consider these
target objects as varying at different scales. We then use the changes of scale to introduce a notion of proximity, in the
form of the associated interleaving distance, \cite{Bub1}, \cite{Bub2}. This is a way of measuring proximity between 
$\cV_\alpha$-type resources by checking whether resource conversion can be inverted after a sufficiently small change
of scale. This provides a convenient measurement of approximation and convergence to the Pareto frontier, with
respect to which one can evaluate possible approximation algorithms. In particular, we will discuss analogs in this
setting of the swarm intelligence algorithms for multi-objective programming. 

\smallskip
\subsection{Scale structure}\label{scaleSec}

As described above, we consider an additional {\em scale} parameter, that we incorporate in the target categories,
by replacing the $\cV_\alpha$ with categories of functors $\cV_\alpha^{(\R,\leq)}={\rm Func}((\R,\leq), \cV_\alpha)$ from
the thin category $(\R,\leq)$ to the category $\cV_\alpha$. Thus, objects in $\cV_\alpha^{(\R,\leq)}$ are determined by a family
$Y_\alpha(s)$ of objects in $\cV_\alpha$, parameterized by $s\in \R$, and a family of morphisms $\varphi_{s\leq s'}: Y(s)\to Y(s')$.
Morphisms $\Hom_{\cV_\alpha^{(\R,\leq)}}(Y,Y')$ are natural transformations, determined by collections of morphisms $\varphi(s): Y(s)\to Y'(s)$ 
in $\cV_\alpha$ satisfying the compatibility $\varphi(s') \varphi_{s\leq s'}=\varphi_{s\leq s'} \varphi(s)$. 

\smallskip

We consider, as before, valuation functors $F_\alpha: \Sigma_\cC(S) \to \cV_\alpha^{(\R,\leq)}$ and target objects $X_\alpha=((X_\alpha(s))_{s\in\R}, \varphi_{s\leq s'})$
in $\cV_\alpha^{(\R,\leq)}$. Thus, in this setting the target objects vary with the scale parameter $s\in \R$. Similarly, using the natural identification of  
functors in  ${\rm Func}(\Sigma_\cC(S), {\rm Func}((\R,\leq), \cV_\alpha))$ with functors in ${\rm Func}(\Sigma_\cC(S)\times (\R,\leq), \cV_\alpha)$, we can
think of the valuation functors themselves as dependent on a scale factor $F_\alpha=(F_{\alpha,s})_{s\in \R}$ with $F_{\alpha,s}:\Sigma_\cC(S) \to \cV_\alpha$.

\smallskip

On the category $\cV_\alpha^{(\R,\leq)}$ there are ``change of scale" functors, which are a special case of the more general flow/coflow structure we
recall in \S \ref{coflowSec} below. For each $\epsilon\geq 0$, there are endofunctors $T_\epsilon: \cV_\alpha^{(\R,\leq)}\to \cV_\alpha^{(\R,\leq)}$ with
$T_0={\rm id}$ and $T_\epsilon T_{\epsilon'}=T_{\epsilon+\epsilon'}$ with 
$T_\epsilon X(s)=X(s+\epsilon)$ and $T_\epsilon\varphi_{s\leq s'}=\varphi_{s+\epsilon\leq s'+\epsilon}$. 

\smallskip

Note that colimits of functors are evaluated pointwise, so if colimits exist in $\cV_\alpha$ then they also exist in $\cV_\alpha^{(\R,\leq)}$. 
The functors $T_\epsilon$ preserve colimits, namely the colimit of the image diagram is the image of the colimit. 

\smallskip
\subsection{Categories with coflow}\label{coflowSec}

We consider here the case where the target categories of the valuation functors are {\em categories with coflow} in the sense of \cite{Cruz}, see
also \cite{deSilva}, \cite{Kashiwara}. Examples of categories with flows (and dually coflows) include persistence modules and
derived sheaves (\cite{deSilva}, \cite{Kashiwara}). We recall briefly the main properties of categories with flows and coflows that we need 
to use in the following.

\smallskip

A {\em flow on a category} $\cV$ is a functor $T: [0,\infty) \to \cE(\cV)$ from the thin category $([0,\infty),\leq)$ to the category $\cE(\cV)$ of
endofunctors of $\cV$, with natural transformations $\mu_{\epsilon, \epsilon'}:T_\epsilon T_{\epsilon'}\to T_{\epsilon +\epsilon'}$ and $\mu_0: {\rm id}_{\cC}\to T_0$,
so that $\mu_{0,\epsilon} \mu_0 {\rm id}_{T_\epsilon} : T_\epsilon \to T_0 T_\epsilon \to T_\epsilon$ and $\mu_{0,\epsilon} {\rm id}_{T_\epsilon} \mu_0  : T_\epsilon \to T_\epsilon T_0  \to T_\epsilon$
are the identity and $\mu_{\epsilon,\zeta+\delta} {\rm id}_{T_\epsilon} \mu_{\zeta,\delta}=\mu_{\epsilon+\zeta,\delta} \mu_{\epsilon,\zeta} {\rm Id}_{T_\delta}$ and
$T_{\epsilon+\zeta\leq \delta+\kappa} \mu_{\epsilon,\zeta}=\mu_{\delta,\kappa} T_{\epsilon\leq \delta} T_{\zeta\leq \kappa}$, with $T_{s\leq s'}$ the natural transformation
in $\cE(\cV)$ associated to the morphism $s\leq s'$ in $[0,\infty)$, see Definition~1 of \cite{Cruz}. 

\smallskip

A {\em coflow on a category} is defined dually. Namely $(\cV,T)$ is a category with coflow iff $(\cV^{op}, T^{op})$ is a category with flow. (The opposite functor $T^{op}$
acts as $T$ on objects and morphisms, but natural transformations between opposite functors have the reversed direction.) The special case of the ``change of scale"
functors on a category $\cV^{(\R,\leq)}$ described in \S \ref{scaleSec} are an example of a {\em strict} flow, where $T_0={\rm id}$ and $T_\epsilon T_{\epsilon'}=T_{\epsilon+\epsilon'}$
instead of having natural transformations between them. While we will work here mostly with this special case, we recall here the more general setting, as most of what we
will discuss generalizes easily to more general categories with coflows.

\smallskip

The main advantage of a flow or coflow structure on a category $\cV$ is that it determines on $X={\rm Obj}(\cV)$ an {\em extended-pseudo-metric}.
We assume here that $X$ is a set. By  extended-pseudo-metric we mean a function $d: X\times X \to \R_{\geq 0}\cup\{\infty\}$ with
\begin{enumerate}
\item $d(x,y)\geq 0$ for all $x,y\in X$, with $d(x,y)=0$ if $x=y$;
\item $d(x,y)=d(y,x)$ for all $x,y\in X$;
\item $d(x,z)\leq d(x,y)+d(y,z)$, for all $x,y,z\in X$.
\end{enumerate}
So unlike an actual metric $d$ can take value $\infty$ and can also take value $0$ on some pairs of non-coincident points. 
The extended-pseudo-metric structure determined by a flow or coflow is called the {\em interleaving distance}. It is defined
as follows. 

\smallskip

As in \cite{Cruz}, \cite{deSilva}, we say that
two objects $A,B$ in $X={\rm Obj}(\cV)$ are $\epsilon$-interleaved if there are morphisms $\alpha: A \to T_\epsilon B$ and $\beta: B \to T_\epsilon A$  in $\cV$
with a commutative diagram
$$ \xymatrix{  T_0 A \ar[dd] & A \ar[l] \ar[dr]^{\alpha\qquad\qquad\qquad\qquad\qquad} & B \ar[dl]_{\qquad\qquad\qquad\qquad\qquad\qquad\beta} \ar[r] & T_0 B \ar[dd] \\
& T_\epsilon A \ar[dr]^{T_\epsilon \alpha\qquad\qquad\qquad\qquad\qquad\qquad\qquad} & T_\epsilon B \ar[dl]_{\qquad\qquad\qquad\qquad\qquad\qquad\qquad\qquad T_\epsilon \beta} & \\
T_{2\epsilon} A & T_\epsilon T_\epsilon A \ar[l] & T_\epsilon T_\epsilon B \ar[r] & T_{2\epsilon} B 
} $$
Then define the {\em interleaving distance} as 
$$ d_{(\cV,T)}(A,B):=\inf  \{ \epsilon\geq 0  \,|\, A \text{ and } B \text{ are $\epsilon$-interleaved} \} \, , $$
were the value can be equal to $\infty$ if no $\epsilon$-interleaving occurs for any $\epsilon \in \R_{\geq 0}$.
The distance is zero if $A$ and $B$ are isomorphic.

\smallskip

We also recall the following result of \cite{Cruz}.

\begin{prop}\label{coflowcolim} %{\rm \cite{Cruz}}
Let $(\cV,T)$ be a category with a coflow with the following properties:
\begin{enumerate}
\item all diagrams $D$ of shape $\bullet \to \bullet \to \bullet \to \cdots$ have colimits in $\cV$;
\item for all $\epsilon >0$, the functor $T_\epsilon$ preserves the colimits, namely if
$T_\epsilon {\rm colim} D \simeq {\rm colim} T_\epsilon D$.
\end{enumerate}
Then $({\rm Obj}(\cV),d_{(\cV,T)})$ is metrically complete.
\end{prop}

\smallskip
\subsubsection{Scale and conversion of resources}\label{scaleresSec}

In the categorical theory of resources \cite{CoFrSp16}, \cite{Fr17}, one interprets morphisms as processes of conversion of resources.
Thus, two isomorphic objects represent resources with the property that there is a conversion process from one to the other that is
fully reversible. Based on this idea, and on the interleaving distance recalled above, we can describe proximity of resources through
the existence of a conversion process that becomes reversible at a different scale. 

\smallskip

Namely, if our category of resources is of the form $\cV^{(\R,\leq)}$ (or is more generally a category with a coflow that we can
think of as implementing changes of scale), we say that resources $A(s)$ and $B(s)$ are $\epsilon$-convertible if there are conversion
processes $\alpha_{s,\epsilon}: A(s) \to B(s+\epsilon)$ and $\beta_{s,\epsilon}: B(s) \to A(s+\epsilon)$ that form an $\epsilon$-interleaving
diagram. We say that $A(s)$ and $B(s)$ are $\epsilon$-close if their interleaving distance $d_{(\cV^{(\R,\leq)},T)}(A,B)\leq \epsilon$, 
with respect to the change of scale functors $T$. 

\smallskip

One can think of this condition, for example, in a setting where higher scales correspond to a coarsening of the system
with fewer averaged out variables, as in a renormalization procedure in statistical physics, for example. In such a setting.
one can have a conversion of resources $A(s)\to B(s)$ at a certain scale $s$, where $B(s)$ is not sufficient to fullly 
reconstruct $A(s)$ but sufficient to reconstruct an averaged out version at a larger scale $A(s+\epsilon)$. In a different
setting, one can instead think of a situation where the system at lower scales is run by simpler variables but at higher scales
it involves more complex emergent phenomena, with a conversion at higher scale not being fully reversible, but still
sufficient to reconstruct the simpler systems at lower scales. In this second case $s\mapsto s+\epsilon$ zooms in to lower scales
rather than to higher scales: we should in general think of $s$ as an order of magnitude and some $\lambda^s$ with $\lambda>0$,
either smaller or larger than $1$, as the actual scale parameter.

\smallskip

\begin{defn}\label{econvert}
In a category of resources of the form $\cV^{(\R,\leq)}$, with a scale parameter $s\in \R$, we say that a conversion of
resources given by a morphism $\varphi_s: A(s) \to B(s)$ is $\epsilon$-reversible if there is a $\beta_{s,\epsilon}: B(s) \to A(s+\epsilon)$
such that $\beta_{s,\epsilon} \circ \varphi_s=T_{\epsilon}: A(s)\to A(s+\epsilon)$.
\end{defn}

\smallskip

\begin{rem}\label{reversed}
An $\epsilon$-reversible morphism $\varphi_s: A(s) \to B(s)$ fits into an $\epsilon$-interleaving diagram
obtained from the two commutative diagrams
$$ \xymatrix{
A(s) \ar[r]^{T_\epsilon\circ \varphi_s\quad} 
\ar@/_2pc/[rr]^{T_{2\epsilon}} 
&  B(s+\epsilon) \ar[r]^{\beta_{s,\epsilon}} & 
A(s+2\epsilon) } $$
and
$$ \xymatrix{ 
B(s) \ar[r]^{\beta_{s,\epsilon}} \ar@/_2pc/[rr]^{T_{2\epsilon}} 
& A(s+\epsilon) \ar[r]^{T_\epsilon\circ \varphi_{s+\epsilon}} & 
B(s+2\epsilon)\, .} $$
\end{rem}

\smallskip
\subsection{Colimits and approximation}\label{colimdSec}

Consider as before a diagram $D$ in ${\rm Diagr}_\alpha$ consisting of a sequence of minorizations as in \eqref{infiniteD}.
Now all the objects are in $\cV_\alpha^{(\R,\leq)}$ so we write the scale dependence explicitly 
\begin{equation}\label{infiniteDscale}
 F_{\alpha,s}(\Phi_0) \stackrel{\varphi_{0,s}}{\rightarrow} F_{\alpha,s}(\Phi_1) \stackrel{\varphi_{1,s}}{\rightarrow} \cdots \rightarrow F_{\alpha,s}(\Phi_n)\stackrel{\varphi_{n,s}}{\rightarrow} \cdots  
\end{equation}
with colimit $L_{\alpha,s}(D)$, for $s\in \R$.

\smallskip

\begin{prop}\label{colimscale}
Given a diagram as in \eqref{infiniteDscale}, if for all $\epsilon>0$ there is an $n_0\in \N$ such that, for all $n\geq n_0$ there 
are maps $\upsilon_{n,\epsilon}: F_{\alpha,s}(\Phi_{n+1}) \to  F_{\alpha,s+\epsilon}(\Phi_n)$ that fit into a commutative diagram
$$ \xymatrix{ F_{\alpha,s}(\Phi_n) \ar[d]_{T_\epsilon} \ar[r]^{\varphi_{n,s}} & F_{\alpha,s}(\Phi_{n+1}) \ar[d]^{T_\epsilon} \ar[dl]^{\upsilon_{n,\epsilon}} \\
F_{\alpha,s+\epsilon}(\Phi_n) \ar[r]_{\varphi_{n,s+\epsilon}} & F_{\alpha,s+\epsilon}(\Phi_{n+1}) } $$
then the sequence $\{ F_{\alpha,s}(\Phi_k) \}$ converges in the interleaving distance to the colimit $L_{\alpha,s}(D)$ of the diagram \eqref{infiniteDscale}.
\end{prop}

\proof We obtain from the commutative diagram above an $\epsilon$-interleaving diagram as in Remark~\ref{reversed}, with
$$ \xymatrix{ F_{\alpha,s}(\Phi_n) \ar[r]^{T_\epsilon\circ \varphi_{n,s}\quad} 
\ar@/_2pc/[rr]^{T_{2\epsilon}} 
&  F_{\alpha,s+\epsilon}(\Phi_{n+1}) \ar[r]^{\upsilon_{n,\epsilon}} & 
F_{\alpha,s+2\epsilon}(\Phi_n) } $$
and
$$ \xymatrix{  F_{\alpha,s}(\Phi_{n+1}) \ar[r]^{\upsilon_{n,\epsilon}} \ar@/_2pc/[rr]^{T_{2\epsilon}} 
& F_{\alpha,s+\epsilon}(\Phi_n)\,\, \,\, \ar[r]^{\varphi_{n,s+2\epsilon}\circ T_\epsilon\quad} & F_{\alpha,s+2\epsilon}(\Phi_{n+1}) } $$
We then have the following commutative diagram, where $\psi_{n,s}$ are the morphisms to the colimit, 
and the morphisms $\omega_{\alpha,\epsilon}: L_{\alpha,s}(D)\to F_{\alpha, s+\epsilon}(\Phi_n)$
are uniquely determined by the universal property of the colimit,
$$ \xymatrix{ F_{\alpha,s}(\Phi_n) \ar[dd]_{T_\epsilon} \ar[r]^{\varphi_{n,s}} \ar[drr]^{\quad\quad\psi_{n,s}} & F_{\alpha,s}(\Phi_{n+1}) \ar[dd]^{T_\epsilon} \ar[ddl]^{\upsilon_{n,\epsilon}} \ar[dr]^{\psi_{n+1,s}} &  \\
& &  L_{\alpha,s}(D) \ar[lld]_{\omega_{\alpha,\epsilon}\quad} \\
F_{\alpha,s+\epsilon}(\Phi_n) \ar[r]_{\varphi_{n,s+\epsilon}} & F_{\alpha,s+\epsilon}(\Phi_{n+1}) & 
 } $$
This similarly determines $\epsilon$-interleaving diagrams
$$ \xymatrix{ F_{\alpha,s}(\Phi_n) \ar[r]^{T_\epsilon\circ \psi_{n,s}\quad} 
\ar@/_2pc/[rr]^{T_{2\epsilon}} 
&  L_{\alpha,s+\epsilon}(D) \ar[r]^{\omega_{\alpha,\epsilon}} & 
F_{\alpha,s+2\epsilon}(\Phi_n) } $$
and
$$ \xymatrix{  L_{\alpha,s}(D) \ar[r]^{\omega_{\alpha,\epsilon}} \ar@/_2pc/[rr]^{T_{2\epsilon}} 
& F_{\alpha,s+\epsilon}(\Phi_n)\,\, \,\, \, \ar[r]^{\psi_{n,s+2\epsilon}\circ T_\epsilon\quad} & L_{\alpha,s+2\epsilon}(D) } $$
Thus, we obtain that the interleaving distance $d(F_{\alpha}(\Phi_n), L_{\alpha}(D))\leq \epsilon$. 
\endproof

\smallskip
\subsection{Particle swarm algorithm}\label{swarm2Sec}

We can then consider a possible particle swarm algorithm designed in the following way.
We are assuming that, given two objects $Y,Y'$ in a target category $\cV_\alpha^{(\R,\leq)}$,
the morphisms $\Hom_{\cV_\alpha^{(\R,\leq)}}(Y,Y')$ are explicitly known. The algorithm also
requires the existence diagrams in ${\rm Diagr}_\alpha$ with $\epsilon$-reversible strict minorizations
for sufficiently large $n$. 
We can then map an $\epsilon$-neighborhood of the Pareto frontier with the following swarm
algorithm, using a large number $n$ of draws for each particle and a 
large number $N$ of particles.

\begin{itemize}
\item Choose an approximation level $\epsilon>0$.
\item Initialize a swarm of $N$ probabilistic particles, by random draws of their initial positions $\Phi_0^{(1)}, \ldots, \Phi_0^{(N)}$
in $\Sigma^{adm}_{\cC}(S)$. 
\item For each $i=1,\ldots, N$ proceed as in the case of the single particle of \S \ref{probpartSec} through successive draws of
positions $\Phi_k^{(i)}$ in $\Sigma^{adm}_{\cC}(S)$, for $k=1,\ldots,n$. 
\item At each successive draw $\Phi_k^{(i)}$ check if the $(F,X)$-strict minorization condition $\Hom_{\cV_\alpha^{(\R,\leq)}}(F_\alpha(\Phi^{(i)}_\ell), F_\alpha(\Phi^{(i)}_k))\neq \emptyset$
holds for previous draws $\Phi_\ell^{(i)}$, $0\leq \ell < k$.
\item If the strict minorization condition holds, check if the $\epsilon$-reversibility, namely the existence of morphisms as in Proposition~\ref{colimscale}, reversing the
minorization direction up to a scale shift.
\item If a strict minorization with $\epsilon$-reversibility exists, then the new draw $\Phi^{(i)}_k$ is in an $\epsilon$-neighborhood of the Pareto frontier.
The chain of minorizations from $F_\alpha(\Phi_0^{(i)})$ to this $F_\alpha(\Phi^{(i)}_k)$ gives a corresponding $\epsilon$-approximation to a diagram in ${\rm Diagr}_\alpha$.
\item If $\epsilon$-invertibility is not satisfied, select a longest chains of strict minorizations in the sequence $\{ F_\alpha(\Phi_0^{(i)}), \ldots, F_\alpha(\Phi_n^{(i)}) \}$. Let 
$F_\alpha(\Phi^{(i)}_k)$ be the last term of this sequence. 
\item Search for chains of strict minorizations starting at $F_\alpha(\Phi^{(i)}_k)$ among the set $\{ F_\alpha(\Phi^{(j)}_k) \}_{j=1}^N$ of the other swarm particles positions at the same time.
\item For each strict minorization check $\epsilon$-invertibility. Whenever an $\epsilon$-reversible strict minorization is found the corresponding particle position is 
in an $\epsilon$-neighborhood of the Pareto frontier.
\item For all the chains of strict minorizations that do not satisfy $\epsilon$-reversibility among the $F_\alpha(\Phi^{(j)}_k)$ continue the search with the new
draws $F_\alpha(\Phi^{(j)}_{k+1})$ and repeat the process.
\end{itemize}

This procedure alternates between new draws $\Phi^{(j)}_{k+1}$ for each particle and searching for best positions for a given particle up to a given time $n$ 
and comparing positions of different particles $\Phi^{(j)}_k$, $i=1,\ldots,N$ at a fixed time $k$, as in the case of the original swarm intelligence algorithm. 
Multiple runs of the algorithm will identify points in an $\epsilon$-neighborhood of the Pareto frontier.

\bigskip
\bigskip
\subsection*{Acknowledgment} The author is partially supported by
NSF grants DMS-1707882 and DMS-2104330
and by FQXi grants FQXi-RFP-1 804 and FQXi-RFP-CPW-2014, 
SVCF grant 2020-224047.

\end{document}